\documentclass{article}
\usepackage{amsfonts,amsmath}
\mathchardef\mhyphen="2D 
\usepackage[a4paper]{geometry}
\usepackage{microtype}
\usepackage{relsize}

\def\zo/{$0\mkern2mu\mhyphen1$}

\def\nn/{$n \times n$}

\DeclareMathOperator\perm{perm}
\title{On the Approximate Asymptotic Statistical Independence of the Permanents of 0-1 Matrices}
\date{\today} 
\author{Paul Federbush\\
Department of Mathematics\\
University of Michigan\\
Ann Arbor, MI, 48109-1043}
\begin{document}

\maketitle
\begin{abstract}
We consider the ensemble of $n \times n$ \zo/ matrices with all column and row sums equal $r$. We give this ensemble the uniform weighting to construct a measure $E$. We know from work of Wanless and Pernici that 
\begin{equation}
\label{A1}
\tag{A1}
E\left(\prod_{i=1}^N \left(\perm_{m_i}(A)\right)\right) = \prod_{i=1}^N \left(E\left(\perm_{m_i}(A)\right)\right) \left(1+\mathcal{O}\left(1/n^4\right)\right) 
\end{equation}
In this paper we prove 
\begin{equation}
\label{A2}
\tag{A2}
E_1\left(\prod_{i=1}^N \left(\perm_{m_i}(A)\right)\right) = \prod_{i=1}^N \left(E_1\left(\perm_{m_i}(A)\right)\right)\left(1+\mathcal{O}\left(1/n^2\right)\right)
\end{equation}
where $E_1$ is the measure constructed on the ensemble of $n \times n$ matrices with non-negative integer entries realized as the sum of $r$ random permutation matrices. $E_1$ is often used as an ``approximation'' to $E$. We have computer evidence for
\begin{equation}
\label{A3}
\tag{A3}
E_1\left(\prod_{i=1}^N \left(\perm_{m_i}(A)\right)\right) = \prod_{i=1}^N \left(E_1\left(\perm_{m_i}(A)\right)\right) \left(1+\mathcal{O}\left(1/n^4\right)\right) 
\end{equation}

\end{abstract}

\section{Introduction}

We consider the ensemble of $n\times n$ \zo/ matrices whose row and column sums all equal $r$. We define the uniform measure in this ensemble, calling it $E$. We know from the work of Wanless \cite{w} and Pernici \cite{p} that
\begin{equation}
\label{1.1}
\tag{1.1}
E\left(\prod_{i=1}^N \left(\perm_{m_i}(A)\right)\right) = \prod_{i=1}^N \left(E\left(\perm_{m_i}(A)\right)\right)\left(1+\mathcal{O}\left(1/n^4\right)\right)
\end{equation}
We let $E_1$ be the measure on $n\times n$ matrices with non-negative integer entries, constructed as the uniform measure on a
sum of $r$ random permutations of $n$ objects. We here prove
\begin{equation}
\label{1.2}
\tag{1.2}
E_1\left(\prod_{i=1}^N \left(\perm_{m_i}(A)\right)\right) = \prod_{i=1}^N \left(E_1\left(\perm_{m_i}(A)\right)\right) \left(1+\mathcal{O}\left(1/n^2\right)\right) 
\end{equation}
In Section 8 we present algebraic ( rigorous ) computer computations for some cases with $N=2, r=2$ and $N=2,r=3$ supporting
\begin{equation}
\label{1.3}
\tag{1.3}
E_1\left(\prod_{i=1}^N \left(\perm_{m_i}(A)\right)\right) = \prod_{i=1}^N \left(E_1\left(\perm_{m_i}(A)\right)\right) \left(1+\mathcal{O}\left(1/n^4\right)\right) 
\end{equation}
These same calculations show (1.3) is not true with $\mathcal{O}\left(1/n^4\right)$ replaced by $\mathcal{O}\left(1/n^5\right)$. 

$E_1$ is often used as an `approximation' to $E$. 
In \cite{f} and \cite{fa} certain aesthetic relations seem to hold in the same form for $E$ and $E_1$ expectations. Similarly we believe eq (\ref{1.1}) presages the truth of eq (1.3). However Section 7 has some results that may cause one some hesitation. 

We briefly consider the Bernoulli random matrix ensemble where each entry independently has a probability $p=r/n$ of being one, and is zero otherwise. We let $E_B$ be the associated measure. $E_B$ is a less worthy "approximation" to $E$, and following along the
lines of the calculation of this paper it is not difficult to show if one replaces $E_1$ in eq(\ref{1.2}) by $E_B$ the resulting relation does not hold.

In a future paper I plan to relate eq (\ref{1.1}) to graph positivity, \cite{bfp}. In particular we plan to use these equations to prove a weak form of the positivity of $\Delta^k d(i)$, see Section II of \cite{bfp}. It was a study of such `graph positivity' that got me involved with the conjecture of this paper. 

One can only appreciate the magic of eq (\ref{1.1}), or eq (\ref{1.2}) , by seeing the complicated calculations and cancellation involved in the proofs. We do not have any understanding 
why the necessary cancellations take place! 

The reader will be forced to embed himself or herself in the world of \cite{f} to follow developments within. But better yet, find another way to attack the study of this conjecture! 

A final note before plunging into the calculation, is the observation that both the conjecture of this paper, and the computation in \cite{f}, can be viewed as an asymptotic statistical independence of the permanents of \zo/ matrices. 

\section{The Strategy}

We depend on the reader being familiar with the first five pages of \cite{f}. In that paper one studies a product of two permanents, here we deal with a product of $N$ permanents, but the ideas are the same. A single permanent $\perm_m(A)$, we view as a sum of `terms', $\binom{n}{m}^2 m!$ such, each term a product of $m$ `entries', $\prod_{i=1}^N \perm_{m_i}(A)$ we view as a sum of `multiterms', each multiterm a product of $N$ `subterms'. Given a fixed multiterm, let $1 \leq i < j \leq N$, and look at the subterms $i$ and $j$. Each entry in the $j$ subterm is in class $1, 2, 3, 4$, or $5$ with respect to the $i$ subterm, in the language of \cite{f}. In this paper we study those types of multiterms in the expansion of $E_1 \left(\prod_{i=1}^N \perm_{m_i}(A)\right)$ that contribute a value larger than $\mathcal{O}\left(1/n^2\right) E_1 \left(\prod_{i=1}^N \perm_{m_i}(A)\right)$. 
More exactly, the sum of all the multiterms of the excluded types is bounded by 
\noindent$\mathcal{O}\left(1/n^2\right) E_1 \left(\prod_{i=1}^N \perm_{m_i}(A)\right)$. The result is that we need consider multiterms that have at most one pair $i < j$ of subterms with one entry of subterm $j$ that is in class other than class 1 to the subterm $i$. In other words, of all the $\sum_{j=1}^N m_j$ entries in the multiterm we consider those that have at most two that share a row or a column! In the following sections we will treat each of the types of multiterms that might contribute to a $1/n$ correction in eq (\ref{1.2}): only class 1 entries, a single class 2 entry, a single class 3 or class 4 entry. Loosely speaking, as is easy to believe, each time two entries share a row or a column one loses a power of $n$. 

\section{Pure Class 1 Multiterms}

We first consider the multiterms where no two entries share a row or a column. In the factor $\left(1+\mathcal{O}\left(1/n^2\right)\right)$ in eq (\ref{1.2}), these multiterms give rise to the $1$, as well as $\mathcal{O}\left(1/n\right)$ corrections canceled to order $1/n^2$ by contributions of the other type multiterms. 

We write down the exact expression for $\prod_{i=1}^N E_1 \left(\perm_{m_i}(A)\right)$, which we call $I$.
\begin{equation}
\label{3.1}
\tag{3.1}
I = \prod_{k=1}^N \left[\underset{\underset{\sum_i m_{k,i}=m_k}{m_{k,i}}}{\sum} \frac{1}{\left(n!\right)^r} \binom{n}{m_k}^2 m_k! \frac{m_k!}{m_{k,1}! \cdots m_{k,r}!} \prod_i \left(n-m_{k,i}\right)!\right]
\end{equation}
$II$ is then the contribution of all pure class 1 multiterms to $E_1\left(\prod_{i=1}^N \left(\perm_{m_i}(A)\right)\right) $. 
\begin{equation}
\label{3.2}
\tag{3.2}
II = \frac{1}{\left(n!\right)^r} \prod_{k=1}^N \left[\underset{\underset{\sum_i m_{k,i}=m_k}{m_{k,i}}}{\sum} \binom{n-\sum_{t<k}m_t}{m_k}^2 m_k! \frac{m_k!}{m_{k,1}! \cdots m_{k,r}!}\right] \prod_i \left(n-\sum_t m_{t,i}\right)!
\end{equation}
We now take all the sums out of the terms $I$ and $II$ writing \begin{equation}
\label{3.3}
\tag{3.3}
I = \sum k~I
\end{equation}
\begin{equation}
\label{3.4}
\tag{3.4}
II = \sum k~ II
\end{equation}
where $k ~ I$ and $k ~ II$ are the `kernels' of the sums. We then write 
\begin{equation}
\label{3.5}
\tag{3.5}
II - I = \sum\left(k~ II - k~I\right) = \sum k~I \left(\frac{k~II}{k~ I}-1\right)
\end{equation}
where 
\begin{align}
\frac{k~II}{k~I}& = \left(n!\right)^{r(N-1)} \frac{\prod_{k=1}^N\binom{n-\sum_{t<k}m_t}{m_k}^2 }{\prod_{k=1}^N \binom{n}{m_k}^2 }\cdot \frac{\prod_{i=1}^r\left(n-\sum_{t}m_{t,i}\right)! }{\prod_{k=1}^N \prod_{i=1}^r \left(n-m_{k,i}\right)! } \label{3.6}\tag{3.6}\\
& = e^\alpha \label{3.7}\tag{3.7}
\end{align}

To evaluate $\alpha$ we turn to the lemma in the appendix, noting that to calculate the $1/n$ terms we need only keep the first term on the right side of eq (A7). This yields up to $\mathcal{O}\left(1/n^2\right)$. 
\begin{align}
\alpha \cong &\frac{1}{2n} \left[2\sum_k \left(\sum_{t<k} m_t\right)^2 + 2\sum_{k} \left(m_k\right)^2-2\sum_k \left(\sum_{t \leq k} m_t\right)^2 \right.\notag\\
& \quad + \left.\sum_i \left(\sum_t m_{t,i}\right)^2-\sum_k \sum_i \left(m_{k,i}\right)^2 \right] \label{3.8}\tag{3.8}
\end{align}
\noindent or, with some thought 
\begin{equation}
\label{3.9}
\tag{3.9}
\alpha \cong \frac{1}{2n} \left[-4\sum_k \sum_{t <k} m_t m_k+2\sum_i \sum_k \sum_{t<k}m_{t,i}m_{k,i} \right]
\end{equation}
In general we write $a \cong b$ to mean $a = b\left(1+\mathcal{O}\left(1/n\right)\right)$. Substituting back into eq (\ref{3.5}) 
\begin{equation}
\label{3.10}
\tag{3.10}
II - I \cong \sum ( k~ I )~\alpha
\end{equation}
To compute the right side of eq (\ref{3.10}) up to $\frac{1}{n}\cdot I$ we can approximate $k~I$ as follows 
\begin{equation}
\label{3.11}
\tag{3.11}
k~I \cong \prod_k \left[\frac{1}{m_{k,1}! \cdots m_{k,r}!} n^{m_k} \right]
\end{equation}
Since 
\begin{align}
& \lim_{n \to \infty} \prod_{k=1}^N \left[\frac{1}{\left(n!\right)^r} \frac{\left(n!\right)^2}{\left(\left(n-m_k\right)!\right)^2} \prod_i \left(n-m_{k,i}\right)! \frac{1}{n^{m_k}}\right]\label{3.12}\tag{3.12}\\
=& \lim_{n \to \infty} \prod_{k=1}^N \left\lbrace\left[\prod_i \left(\frac{\left(n-m_{k,i}\right)!}{n!}\right)n^{m_k} \right]\cdot \left[\frac{1}{n^{m_k}}\left(\frac{n!}{(n-m_k)!}\right)\right]^2 \right\rbrace \label{3.13}\tag{3.13}\\
=&~ 1 \label{3.14}\tag{3.14}
\end{align}
\noindent where in fact each of the expressions in brackets in eq (\ref{3.13}) approaches $1$ with $n$ going to infinity. So finally 
\begin{equation}
\label{3.15}
\tag{3.15}
II -I \cong \sum\left(\prod_{k=1}^N \left(\frac{n^{m_k}}{m_{k,1}!\cdots m_{k,r}!}\right) \right) \frac{1}{2n} \left[-4 \sum_k \sum_{t <k}m_t m_k + 2\sum_i \sum_k \sum_{t<k} m_{t,i}m_{k,i} \right]
\end{equation}
We now use eq (A2) from the appendix to arrive at our final result 
\begin{equation}
\label{3.16}
\tag{3.16}
II-I \cong \left[\prod_{k=1}^N\left(\frac{n^{m_k}r^{m_k}}{m_k!}\right) \right] \cdot \frac{1}{2n} \left[-4\sum_k \sum_{t <k}m_t m_k + 2\sum_i \sum_k \sum_{t<k} \frac{m_t}{r} \frac{m_k}{r} \right]
\end{equation}
or 
\begin{equation}
\label{3.17}
\tag{3.17}
II-I = \left[\prod_{k=1}^N  \left(\frac{n^{m_k} r^{m_k}}{m_k!}\right) \right]\frac{1}{2n} \left[-4\sum_k \sum_{t<k}\left(m_t m_k\right)\left(1-\frac{1}{2r}\right) \right] + I\mathcal{O}\left(1/n^2\right)
\end{equation}

Alternatively we could avoid the discussion from eq (\ref{3.11}) to eq (\ref{3.17}), going directly from (\ref{3.10}) to 
\begin{equation}
\label{3.18}
\tag{3.18}
II-I = I \cdot \frac{1}{2n} \left[-4\sum_k \sum_{t<k} \left(m_t m_k\right)\left(1-\frac{1}{2r}\right) \right] + I\mathcal{O}\left(1/n^2\right)
\end{equation}
using 
\begin{equation}
\label{3.19}
\tag{3.19}
\underset{\underset{\sum m_{k,i} = m_k}{m_{k,i}}}{\sum} \frac{1}{\left(n!\right)^r} \binom{n}{m_k}^2 m_k! \frac{m_k!}{m_{k,1}! \cdots m_{k,r}!} \prod_i \left(n-m_{k,i}\right) \cdot \left(m_{k,j}-\frac{m_k}{r}\right) =0
\end{equation}
that follows from the symmetry of the $r$ colors. 

\section{A Single Entry of Class 2}

We now study the contribution to $E_1 \left(\prod_{i=1}^N \perm_{m_i}(A)\right)$ of multiterms with all entries but one of class 1 and the single other entry of class 2. We sum over the color $s$ of this entry, the subterm $j$ it is in, and the subterm $i$ it is in class 2 relative to. We let this contribution be called $III$. \begin{align}
III = &\sum_s \sum_{i<j} \frac{1}{\left(n!\right)^r} \prod_{k=1}^{j-1} \left[\underset{\underset{\sum m_{k,i} = m_k}{m_{k,i}}}{\sum} \binom{n-\sum_{t<k}m_t}{m_k} ^2m_k! \frac{m_k!}{m_{k,1}! \cdots m_{k,r}!} \right]\cdot \notag\\
& \quad \cdot \left[\underset{\underset{\sum m_{j,i} = m_j-1}{m_{j,i}}}{\sum}\binom{n-\sum_{t<j}m_t}{m_j-1}^2 \left(m_j -1\right)! \frac{\left(m_j-1\right)!}{m_{j,1}! \cdots m_{j,r}!} \right] \cdot \notag\\
& \prod_{k=j+1}^N \left[\underset{\underset{\sum m_{k,i} = m_k}{m_{k,i}}}{\sum}\binom{n-\sum_{t<k}m_t+1}{m_k}^2m_k!\frac{m_k!}{m_{k,1}! \cdots m_{k,r}!} \right]\cdot \notag\\
& \prod_i (n-\sum_t m_{t,i})
!\, m_{i,s} \label{4.1}\tag{4.1}\\
\end{align}
The final $m_{i,s}$ represents the choice of which entry of color $s$ in subterm $i$ agrees with the entry in subterm $j$. In subterm $j$ there are just $m_j -1$ entries to be summed over once the class 2 entry has been selected. We write (\ref{4.1}) as (we will use $Q_j$ later) 
\begin{equation}
=  \sum_s \sum_{i<j} Q_j \prod_i \left(n-\sum_tm_{t,i}\right)m_{i,s} \label{4.2}\tag{4.2}
\end{equation}
 and here also define $II_j'$ 
\begin{equation}
\label{4.3}
\tag{4.3}
III = \sum_s \sum_{i<j} II_j' \, m_{i,s}
\end{equation}
where $II_j'$ is exactly $II$ of eq (\ref{3.2}) with $m_j$ replaced by $m_j-1$. We now use the argument as in eq (\ref{3.11})-(\ref{3.14}) to write 
\begin{align}
II_j' \cong & \prod_{\underset{k \neq j}{k=1}}^N \left(\underset{\underset{\sum m_{k,i} = m_k}{m_{k,i}}}{\sum} \frac{1}{m_{k,1}! \cdots m_{k,r}!} n^{m_k}\right) \cdot \underset{\underset{\sum m_{j,i} = m_j-1}{m_{j,i}}}{\sum} \left(\frac{1}{m_{j,1}! \cdots m_{j,r}!} n^{m_j-1}\right) \label{4.4}\tag{4.4}\\
\cong & I \cdot \frac{m_j}{nr} \label{4.5}\tag{4.5}
\end{align}
by eq (A.1) from the appendix. Substituting into eq (\ref{4.3}) we get 
\begin{align}
III \cong &\sum_s \sum_{i<j} I \frac{m_j}{nr}m_{i,s}\notag\\
\cong & I~ \sum_{i<j} \frac{m_j m_i}{nr}\label{4.6}\tag{4.6}
\end{align}
equivalently
\begin{equation}
\label{4.7}
\tag{4.7}
III = I \, \left(\sum_{i<j} \frac{m_j m_i}{nr} + \mathcal{O}\left(1/n^2\right)\right)
\end{equation}

\section{A Single Class 3 or Class 4 Entry}

We calculate the contribution of the sum of these two types of multiterms, multiplying by 2 the contribution of the multiterms with a single class 3 entry. We let $IV$ be the sum of these two types. We specify all the class 1 entries in the multiterm and at the last choice select the single class 3 entry, arriving at 
\begin{equation}
\label{5.1}
\tag{5.1}
IV = 2\sum_{a,b} \sum_{i<j} Q_j m_{i,a} \left(n-\sum_t m_t +1\right) \prod_i (n-\sum_t m_{t,i})! 
\end{equation}
see eq (\ref{4.2})

The 2 is for the sum of the two classes, $b$ is the color of the class 3 or class 4 entry, $a$ is the color of the entry in subterm $i$ that shares a row or column with this entry. $m_{i,a}$ selects the color $a$ entry in subterm $i$ sharing the row or column, $(n-\sum_t m_t+1)$ selects the column (row) for the class 3 (class 4) entry. 
\begin{equation}
\label{5.2}
\tag{5.2}
IV \cong 2 (r-1) \sum_{i<j} I \, m_i \frac{m_j}{nr}
\end{equation}
$(r-1)$ arises choosing the color of the class 3 or 4 entry. Finally 
\begin{equation}
\label{5.3}
\tag{5.3}
IV = 2 (r-1) \sum_{i<j} I \, \frac{m_i m_j}{rn} \left(1+\mathcal{O}\left(1/n^2\right)\right)
\end{equation}

\section{Summing the Terms}

From eq (\ref{3.17}) we get 
\begin{equation}
\label{6.1}
\tag{6.1}
II - I\cong \sum_{i<j} I \left(-\frac{2}{n} \left(1-\frac{1}{2r}\right)\right) m_i m_j
\end{equation}
From eq (\ref{4.6}) we have
\begin{equation}
\label{6.2}
\tag{6.2}
III \cong \sum_{i<j} I \, \left(\frac{m_i m_j}{nr}\right) 
\end{equation}
and from eq (\ref{5.3}) we have 
\begin{equation}
\label{6.3}
\tag{6.3}
IV \cong 2(r-1) \sum_{i<j} \left(\frac{m_i m_j}{rn}\right) I
\end{equation}
Therefore 
\begin{equation}
\label{6.4}
\tag{6.4}
II + III + IV \cong I
\end{equation}

We remark to the reader that a single class 5 entry makes no contribution to order $1/n$. If one looks at the details of the computations of this paper addressed for getting out corrections to order $1/n$, one can see how complex it will be to get the corrections to order $1/n^2$. We have used integer arithmetic computer computations for $N=2, r=2$ and $ r=3$, to help check the correctness of analytic expressions. 

\section{A Cautionary Calculation}

Whereas we have referred to a couple instances where $E_1$ expectations have the same behavior as $E$ expectations, we here present a blatant contrast. In \cite{p}, eq (11) and Appendix A, there is presented an expansion in descending powers of $n$ 
\begin{align}
E\left(\perm_{m}(A)\right)& = an^m + bn^{m-1}+cn^{m-2}+ \cdots \label{7.1}\tag{7.1}\\
a & = (r^m/m!) \label{7.2}\tag{7.2}\\
b & =(r^m/m! )\left(m(m-1)\right) \left(-1+\frac{1}{2}\frac{1}{r}\right)\label{7.3}\tag{7.3}\\
c & = (r^m/m! )\left(m(m-1)(m-2)\right) \left(\frac{1}{6}(3m+1)-\frac{1}{2}(m+1)\frac{1}{r}+\frac{1}{24}(3m+7)\frac{1}{r^2}\right) \label{7.4}\tag{7.4}
\end{align}

(M. Pernici has shown me how to alternatively derive these equations using the formalism of \cite{fa}, assuming Conjecture 1 therein.)

It is easy to check that the corresponding expansion for $E_1 \left(\perm_m(A)\right)$ agrees in the first two terms, but \underline{not} in the third, checked by computer to some high order. We take this expression for  $E_1 \left(\perm_m(A)\right)$ as established. 
This is a clash to the same power of $n$ as pursued in the sequel, paper II. 

\section{Some Computer Calculations}

We note that by the result of this paper (1.2) and the expressions of Section 7 one has
\begin{equation}
\label{8.1}
\tag{8.1}
E_1\left(\prod_{i=1}^N \left(\perm_{m_i}(A)\right)\right) \cong c n^ { \sum {m_i}}
\end{equation}
where $c$ depends on $r$ and the $m_i$. Similarly from (1.1) and from the expressions from Section 7, one has
\begin{equation}
\label{8.2}
\tag{8.2}
E\left(\prod_{i=1}^N \left(\perm_{m_i}(A)\right)\right) \cong c n^ { \sum {m_i}}
\end{equation}
with the value of $c$ the same in (8.1) and (8.2).

We have computed exactly, in the case $r=3$ ( and $N=2$ )

\begin{equation}
\label{8.3}
\tag{8.3}
E_1(perm_{m_1}(A) perm_{m_2}(A)) - E_1(perm_{m_1}(A)) E_1(perm_{m_2}(A))
\end{equation}
for all $m_1,m_2 < 7$ and found that the highest non-vanishing term in the expansion in powers of $n$ was $m_1+m_2-4$. We have done the same for $r=2$ ( and $N=2$ ) for all $m_1,m_2 < 9$ and found that the highest non-vanishing term in the expansion in powers of $n$ was $m_1+m_2-4$. For {r=2} the expectations in (8.3) are polynomials in {n}, for {r=3} rational functions of {n}.

We will study expectations of $perm_5(A)perm_3(A)$ as an example. We are working with an $N=2$ example so that we can compute
the $E_1$ expectation exactly, and with $r=2$.
( It is possible to study expectations in the $E_1$ measure for $N > 2$, but at a tremendous
increase in computational complexity. ) We define
\begin{equation}
\label{8.4}
\tag{8.4}
Q1 \equiv E_1( perm_{5}(A) perm_{3}(A)) 
\end{equation}
\begin{equation}
\label{8.5}
\tag{8.5}
Q2 \equiv E_1( perm_{5}(A) perm_{3}(A)) -  E_1(perm_{5}(A)) E_1(perm_{3}(A))
\end{equation}
For $Q1$ we have an exact expression
\begin{equation}
\label{8.6}
\tag{8.6}
\begin{aligned}
Q1 = &(16/45)n^8  - (104/15)n^7   +(868/15)n^6-(4046/15)n^5+(11548/15)n^4 \\
&-(20926/15)n^3+(73076/45)n^2-(17824/15)n+448 
\end{aligned}
\end{equation}
Similarly we have an exact expression for $Q2$
\begin{equation}
\label{8.7}
\tag{8.7}
Q2=( 8/3  )n^4- ( 100/3 )n^3+(460/3  )n^2-( 920/3 )n +224
\end{equation}

\section*{Appendix}
\begin{equation}
\label{A.1}
\tag{A.1}
\sum_{\underset{\sum m_i = m}{m_i}} \frac{1}{m_1! \cdots m_r!} = \frac{(x_1 + \cdots + x_r)^m}{m!}\bigg|_{x_i =1} = \frac{r^m}{m!}
\end{equation}
\begin{equation}
\label{A.2}
\tag{A.2}
\sum_{\underset{\sum m_i = m}{m_i}} \frac{m_1}{m_1! \cdots m_r!} = x_1 \frac{d}{dx_1}\left( \frac{(x_1 + \cdots + x_r)^m}{m!}\right)\bigg|_{x_i =1} =\frac{m}{r} \frac{r^m}{m!}
\end{equation}
\begin{equation}
\label{A.3}
\tag{A.3}
\sum_{\underset{\sum m_i = m}{m_i}} \frac{m_1^2}{m_1! \cdots m_r!} = x_1 \frac{d}{dx_1}x_1 \frac{d}{dx_1}\left( \frac{(x_1 + \cdots + x_r)^m}{m!}\right)\bigg|_{x_i =1} =m(m-1)\frac{r^{m-2}}{m!}+m\frac{r^{m-1}}{m!}
\end{equation}
\begin{equation}
\label{A.4}
\tag{A.4}
\sum_{\underset{\sum m_i = m}{m_i}} \frac{m_1m_2}{m_1! \cdots m_r!} = x_1 \frac{d}{dx_1}x_2 \frac{d}{dx_2}\left( \frac{(x_1 + \cdots + x_r)^m}{m!}\right)\bigg|_{x_i =1} =m(m-1)\frac{r^{m-2}}{m!}
\end{equation}
\begin{equation}
\label{A.5}
\tag{A.5}
\ln n! = n\ln n -n + \frac{1}{2}\ln(2\pi) + \frac{1}{2}\ln n + \frac{1}{12}\frac{1}{n} + \mathcal{O}\left(\frac{1}{n^3}\right)
\end{equation}

\textbf{\underline{Lemma}}\quad Suppose \begin{equation}
\label{A.6}
\tag{A.6}
\sum a_i = 0 \quad \text{and} \quad \sum a_i q_i =0
\end{equation}
then \begin{equation}
\label{A.7}
\tag{A.7}
\sum a_i \ln \left((n-q_i)!\right) = \sum a_i \left[\frac{1}{2}\frac{q_i^2}{n} + \frac{1}{6}\frac{q_i^3}{n^2}-\frac{1}{4}\frac{q_i^2}{n^2}\right] + \mathcal{O}\left(1/n^3\right)
\end{equation}

The lemma follows from eq(A.5) with a little calculation.

\end{document}